\documentclass[11pt]{article}
\usepackage[latin1]{inputenc}
\usepackage{amsmath}
\usepackage{amsthm}
\usepackage{amsfonts}
\usepackage{amssymb}
\usepackage{graphicx}

\usepackage{amsmath,amsthm,amssymb,graphicx, multicol, array}
\usepackage{enumerate}
\usepackage{enumitem}
\usepackage{hyperref}

\newtheorem{thm}{Theorem}[section]
\newtheorem{lem}{Lemma}[section]

\newtheorem{rem}{Remark}[section]

\newcommand{\N}{\mathbb{N}}

\newcommand{\C}{\mathbb{C}}

\usepackage{fancyhdr}
\begin{document}
\thispagestyle{empty}
\title{\textbf{A Simple Proof Of The Prime Number Theorem} }
\author{N. A. Carella}
\date{}
\maketitle

\textbf{Abstract:} It is shown that the \textit{Mean Value Theorem} for arithmetic functions, and simple properties of the zeta function are sufficient to assemble proofs of the \textit{Prime Number Theorem}, and \textit{Dirichlet Theorem} for primes in arithmetic progressions. These are among the simplest proofs of the asymptotic formulas for the corresponding prime counting functions.
\\

\textbf{Mathematics Subject Classifications: }11A41, 11N05, 11N13.\\
\textbf{Keywords:} Prime Numbers, Distribution of Primes, Prime Number Theorem.

\section{Introduction}
The Prime Number Theorem is an asymptotic formula $\pi(x)=x/\log x +o(x/\log x)$ for the number of primes $p\geq2$ up to a number $x \geq 1$. The quantity $E(x)=o(x/\log(x)$ is called the error term. One of the objectives of prime number theory is to reduce the error term to the optimal $E(x)=O(x^{1/2}\log^2(x)$. Generally, the various proofs of the Prime Number Theorem are derived from the zerofree regions of the zeta function $\zeta(s)$ of a complex number $s\in \C$. The first proofs, independently achieved by Hadamard, and delaVallee Poussin, are based on the zerofree region $\{s\in \C: \mathcal{R}e(s) \geq 1 \}$. Since then, many other related proofs have been found. Furthermore, the zerofree regions of the zeta function, and the error term  have been slightly improved to $E(x)=O(xe^{-c(\log \log x)^{1/3}})$, where $c>0$ is an absolute constant, see \cite{FK02}. The exception to this rule are the more recently discovered elementary proofs of the Prime Number Theorem, which are derived from the average orders of various arithmetic functions. The elementary methods do not require any information on the zeta function, refer to \cite{SA49}, \cite{ES49}, and \cite{DS70}. \\

A simple proof of the Prime Number Theorem is constructed from Mean Value Theorem for arithmetic functions, and basic properties of the zeta function. This proof does not require any deep knowledge of the prime numbers, and it does not require any difficult to prove zerofree region beyond the well known, and simplest zerofree region $\{s\in \C: \mathcal{R}e(s) > 1 \}$. In addition, a similar proof of the Dirichlet Theorem for primes in arithmetic progressions is included.

\begin{thm} \label{thm1.1}   {\normalfont (Prime Number Theorem) }  For all sufficiently large number $x \geq 1$, the prime counting function satisfies the asymptotic formula
\begin{equation}
\pi(x)=\#\{p \leq x: p \;{\normalfont prime} \}=\frac{x}{\log x}+o\left ( \frac{x}{\log x} \right ).
\end{equation}
\end{thm}

\begin{proof} For $n \geq 1$, let $f(n)=\Lambda(n)$ be the weighted characteristic function of the prime numbers, see the definition (\ref{eq220.53}). The corresponding generating function is
\begin{equation}
\sum_{n \geq 1}\frac{\Lambda(n)}{n^s}=- \frac{\zeta^{'}(s)}{\zeta(s)}=-\zeta(s) \cdot \frac{\zeta^{'}(s)}{\zeta^2(s)}
\end{equation}
for $s\in \C$ with $\mathcal{R}e(s)>1$, see Lemma \ref{lem220.01}. By Lemma \ref{lem220.04}, the function 
\begin{equation}
\frac{\zeta^{'}(s)}{\zeta^2(s)}=\sum_{n \geq 1}\frac{g(n)}{n^s}
\end{equation}

converges at $s = 1$. Therefore, by Theorem \ref{thm3399.03} or Theorem \ref{thm3399.09}, the mean value of $f(n)=\Lambda(n)$, in Section \ref{s3399},  is given by
\begin{equation}
M(f)= \lim_{x \to \infty}\frac{1}{x}\sum_{n \leq x}\Lambda(n)=- \frac{\zeta^{'}(1)}{\zeta^2(1)}=1,
\end{equation}

see equations (\ref{eq220.22}) and (\ref{eq220.23}). These immediately imply that
\begin{equation}
\sum_{n \leq x}\Lambda(n)=x+o(x).
\end{equation}
Next, let $\psi(t)=\sum_{n \leq t}\Lambda(n)$. Then, by partial summation, it follows that
\begin{eqnarray}
\pi(x)=\sum_{p\leq x}1&=&\sum_{n \geq x}\frac{\Lambda(n)}{\log n}+O\left( x^{1/2}\log^2 x \right ) \nonumber \\
&=&\int_2^x\frac{1}{\log t} d\psi(t)+O\left( x^{1/2}\log^2 x \right ) \\
&=& \frac{x}{\log x}+o\left( \frac{x}{\log x} \right )\nonumber,
\end{eqnarray}
here the error term $O\left( x^{1/2}\log^2 x \right ) $ accounts for the prime powers $p^k$ with $k \geq 2$.       
\end{proof}

This is probably one of the simplest proof of the Prime Number Theorem, there is a claim for the simplest proof in \cite[ p.\ 80]{DL12}.

\subsection{Extension to Dirichlet Theorem}
As the mean value theorem for arithmetic functions has a natural extension to arithmetic progression, it is natural to consider the next result.

\begin{thm} \label{[thm2.1}     {\normalfont (Dirichlet Theorem) }  For all sufficiently large number $x\geq 1$, and a pair of integers $a\geq 1$ and $q \geq 1$, with $\gcd(a,q)=1$, the prime counting function satisfies the asymptotic formula
\begin{equation}
\pi(x,q,a)=\#\{p \leq x: p \equiv a \bmod q\} \geq \frac{1}{q} \left (1+o\left (1\right ) \right )\frac{x}{\log x}+o\left ( \frac{x}{q\log x} \right ).
\end{equation}
\end{thm}

\begin{proof} Assume $q\geq 2$ is prime, and consider the Mobius pair 
\begin{equation}
f(n)= \sum_{d \mid n} g(d)   \quad  \text{ and } \quad g(n)= -\mu(n) \log n.
\end{equation}

By Theorem \ref{thm4902.31}, the mean value of the arithmetic function $f(n)=\Lambda(n)$ over arithmetic progression $\{qn+a:\gcd(a,q)=1 \text{ and } n\geq 0\}$ is given by
\begin{eqnarray}
M(f)&=& \lim_{x \to \infty}\; \frac{1}{x}\sum_{n \leq x, \; n \equiv a \bmod q}\Lambda(n) \nonumber \\
&=& \frac{1}{q} \sum_{d\mid q} \frac{c_d(a)}{d} \sum_{n \geq 1} \frac{g(dn)}{n} \nonumber \\
&=&- \frac{1}{q} \sum_{d\mid q} \frac{c_d(a)}{d} \sum_{n \geq 1} \frac{\mu(dn) \log (dn) }{n} \nonumber \\
&=&- \frac{c_1(a)}{q} \sum_{n \geq 1} \frac{\mu(n) \log n }{n} - \frac{1}{q} \sum_{1<d\mid q} \frac{c_d(a)}{d} \sum_{n \geq 1} \frac{\mu(dn) \log (dn) }{n} \nonumber \\
&=& \frac{1}{q} \left (- \frac{\zeta^{'}(1)}{\zeta(1)^2} -  \sum_{1<d\mid q} \frac{\mu(d)}{d} \sum_{n \geq 1} \frac{\mu(dn) \log (dn) }{n} \right )\nonumber \\
&=& \frac{1}{q} \left (1- \sum_{1<d\mid q} \frac{\mu(d)}{d} \sum_{n \geq 1} \frac{\mu(dn) \log (dn) }{n} \right ).
\end{eqnarray}
The last line follows from Lemma \ref{lem220.04}, and $c_1(a)=1$. This and Theorem \ref{thm3399.09} immediately imply that
\begin{equation}
\pi(x,q,a)\geq \frac{1}{q} \left ( 1+o\left ( 1\right ) \right )\frac{x}{\log x}+o\left ( \frac{x}{q\log x} \right ),
\end{equation}
where $M(f)=1/q+o(1/q)$. This proves the claim.                                                  
\end{proof}

\begin{rem} {\normalfont  The proof of Theorem \ref{thm4902.31} has no limitations on the range of values of $q \geq 1$. Thus, this result is probably an improvement on the Siegel-Walfisz Theorem, which states that
\begin{equation}
\pi(x,q,a)=\frac{1}{\varphi(q)} \frac{x}{\log x} +O\left (xe^{-c(\log \log x)^{\beta}} \right ),
\end{equation}
where $q=O(\log^B x)$, with $B > 0$, $c>0$ is an absolute constant, and $0 < \beta <1$ constants.}
\end{rem}

\section{Powers Series Expansions of the Zeta Function} \label{s220}
Let $\N=\{ 0, 1, 2, 3, \ldots \}$ be a nonnegative integer. The Mobius function is defined by
\begin{equation}
\mu(n)= 
\begin{cases}
    1       & \quad \text{if } n \text{ is squarefree},\\
    0 & \quad \text{if } n \text{ is not squarefree}.\\
  \end{cases}
\end{equation}

And the vonMangoldt function is defined by
\begin{equation} \label{eq220.53}
\Lambda(n)= 
\begin{cases}
    \log p       & \quad \text{if } n=p^k, k \geq 1, \text{ is a prime power},\\
    0 & \quad \text{if } n\ne p^k, k \geq 1, \text{ is not a prime power}.\\
  \end{cases}
\end{equation}

\begin{lem} \label{lem220.01}   Let $\Lambda$ be the vonMangoldt function, and let $\zeta(s)=\sum_{n \geq 1}n^{-s}$ be the zeta function. Then, 
\begin{equation}
- \frac{\zeta^{'}(s)}{\zeta(s)}=\sum_{n \geq 1}\frac{\Lambda(n)}{n^s}
\end{equation}
is absolutely convergent on the complex half plane $\mathcal{R}e(s) > 1$.   
\end{lem}

\begin{proof} Since the zeta function is convergent on the complex half plane $\{s \in \C:\mathcal{R}e(s) > 1\}$, and it has an Euler product
\begin{equation}
\zeta(s)=\sum_{n \geq 1}\frac{1}{n^s}=\prod_{p \geq 2} \left(1 -\frac{1}{p^s}\right )^{-1}.
\end{equation}
Taking the logarithm derivative yields
\begin{eqnarray}
\frac{d}{ds} \log \zeta(s) &=&-\sum_{p \geq 2} \frac{d}{ds} \log \left(1 -\frac{1}{p^s} \right ) \\
&=&-\sum_{p \geq 2} \frac{d}{ds} \sum_{n  \geq 1}\frac{1}{np^{ns}}  \nonumber \\
&=&\sum_{p \geq 2} \sum_{n  \geq 1}\frac{\log p}{p^{ns}}  \nonumber \\
&=&-\sum_{n  \geq 1}\frac{\Lambda(n)}{n^{s}}  \nonumber .
\end{eqnarray}
The second line follows from the absolute converges on the complex half plane $\mathcal{R}e(s) > 1$, rearranging the double sums, and the definition of the vonMangoldt function in equation (\ref{eq220.53}).  
\end{proof}

The zeta function $\zeta(s)=\sum_{n \geq 1}n^{-s}$ is an entire function of a complex variable $s \in \C$ with a simple pole at $s = 1$. At each fixed complex number $s_0 \in \C$, it has a Taylor series expansion of the form
\begin{equation}
\zeta(s)=-\frac{1}{s-1}+\sum_{n \geq 0}\frac{\zeta^{(n)}(s_0)}{n!}(s-s_0)^{n},
\end{equation}
where $\zeta^{(n)}(s)$ is the $n$th derivative.\\

The power series expansions at the zeros, and some other special values of the zeta function $\zeta(s)$ are simpler to determine, for example, $s_0 \in \mathcal{Z}=\{-2n : n \geq 0 \}$. The well known expansion at $s = 1$ has the power series 
\begin{eqnarray} \label{eq220.03}
\zeta(s)&=&-\frac{1}{s-1}+\sum_{n \geq 0}\frac{\gamma_n}{n!}(s-1)^{n} \\
&=&-\frac{1}{s-1}+\gamma_0-\gamma_1(s-1)+\frac{1}{2} \gamma_1(s-1)^2+\cdots\nonumber,
\end{eqnarray}
see where $\gamma_n$ is the $n$th Stieltjes constant, see \cite{KJ92}, and \cite[Eq. 25.2.4]{DLMF}. And the expansion at $s = 0$ has the power series
\begin{equation}
\zeta(s)=-\frac{1}{s-1}+\sum_{n \geq 0}(-1)^n \delta_n s^{n}, 
\end{equation}                                                
see \cite{DLMF}. The power series expansions at the zeros $s_0$ are some sort of `zeta modular forms' since the first coefficient $a_0 = \zeta(s_0) = 0$.

\begin{lem} \label{lem220.04} The complex value function $-\zeta^{'}(s)/\zeta^2(s)$ is analytic and zerofree on the complex half plane $\Re e(s) \geq 1$, and at $s = 1$. 
\end{lem}

\begin{proof} Since the zeta function is zerofree on the complex half plane $\Re e(s) > 1$, it is sufficient to consider it at the simple pole at $s = 1$. Toward this end, use the power series (\ref{eq220.03}) to rewrite the function as a ratio
\begin{equation}  \label{eq220.22}
- \frac{\zeta^{'}(s)}{\zeta^2(s)}=\frac{-(s-1)^{-2}-\gamma_1+\gamma_2(s-1)+ \cdots}
{(s-1)^{-2}-\gamma_1(s-1)+\gamma_2(s-1)^2+ \cdots}.
\end{equation}
\
Taking the limit yields
\begin{equation}  \label{eq220.23}
 \frac{\zeta^{'}(1)}{\zeta^2(1)}=\lim_{s \to 1}\frac{(s-1)^2\zeta^{'}(s)}{(s-1)^2\zeta^2(s)}=-1.
\end{equation}
This shows that  is well defined at each point on the complex half plane $\Re e(s) > 1$, and at $s = 1$.                                                                         \end{proof}

\begin{lem} \label{lem220.08}   Let $\mu, \Lambda : \N \to \C$ be the Mobius, and the vonMangoldt arithmetic functions. Then
\begin{equation}
\sum_{n \leq x} \sum_{d\mid n} \mu(d) \Lambda(n/d)=o(x).
\end{equation}
\end{lem}

\begin{proof} By Lemma \ref{lem220.01} and  Lemma \ref{lem220.04}, the series
\begin{eqnarray}
- \frac{\zeta^{'}(s)}{\zeta(s)} \cdot \frac{1}{\zeta(s)}&=&-\sum_{n \geq 1}\frac{\Lambda(n)}{n^s} \cdot \sum_{n \geq 1}\frac{\mu(n)}{n^s} \\
&=&-\sum_{n \geq 1}\frac{\sum_{d\mid n} \mu(d) \Lambda(n/d)}{n^s} \nonumber\\
&=&\sum_{n \geq 1}\frac{f(n)}{n^s} \nonumber
\end{eqnarray}
converges on the complex half plane $\Re e(s) > 1$, and at $s = 1$. Therefore, by Lemma \ref{lem3399.01}, the summatory function
\begin{equation}
\sum_{n \leq x}f(n)=\sum_{n \leq x} \sum_{d\mid n} \mu(d) \Lambda(n/d)=o(x).
\end{equation}
for large $x\geq 1$.                                                                            
\end{proof}

\section{Mean Values of Arithmetic Functions} \label{s3399}
Let $f : \N \to \C$ be a complex valued arithmetic function on the set of nonnegative integers. The mean value of an arithmetic function is defined by
\begin{equation} \label{eq3399.003}
M(f)= \lim_{x \to \infty}\,\frac{1}{x}\sum_{n \leq x}f(n).
\end{equation}
The mean value of an arithmetic function is sort of a \textit{weighted density} of the subset of integers $supp(f)=\{n\in \N: f(n)\ne 0\} \subset \N$, which is the \textit{support} of the function $f$, see \cite[p.\ 46]{SS94} for a discussion of the mean value. The \textit{natural density} of a subset of integers $\mathcal{A} \subset \N$ is defined by
\begin{equation}
\delta(\mathcal{A})= \lim_{x \to \infty}\,\frac{\{n \leq x: n \in \mathcal{A}\}}{x}.
\end{equation}

\subsection{Some Results on Mean Values of Arithmetic Functions}
This section investigates the two cases of convergent series, and divergent series
\begin{equation}
\sum_{n \geq 1}\frac{f(n)}{n} < \infty  \qquad \mbox{ and } \qquad  \sum_{n \geq 1}\frac{f(n)}{n} =\pm \infty,
\end{equation}
and the corresponding mean values
\begin{equation}
M(f)= \lim_{x \to \infty} \frac{1}{x}\sum_{n \leq x}f(n) =0  \qquad \mbox{ and } \qquad M(f)= \lim_{x \to \infty} \frac{1}{x}\sum_{n \leq x}f(n) \ne 0 \nonumber
\end{equation}
respectively. First, a result on the case of convergent series is considered here.

\begin{lem} \label{lem3399.01} Let $f : \N \longrightarrow \C$ be an arithmetic function. If the series  $\sum_{n \geq 1}f(n)n^{-1}$ converges, then its mean value  
\begin{equation}
M(f)= \lim_{x \to \infty} \frac{1}{x}\sum_{n \leq x}f(n) =0  
\end{equation}
vanishes.     
\end{lem}
                            
\begin{proof} Consider the pair of finite sums $\sum_{n \leq x}f(n) $, and $\sum_{n \leq x}f(n)n^{-1} $. By hypothesis, $\sum_{n \leq x}f(n)n^{-1}=c+o(1) $, where $c\ne 0$ is constant, for large $ x \geq  1$. Therefore
\begin{eqnarray}
\sum_{n \leq x}f(n) &=& \sum_{n \geq 1}n \cdot \frac{f(n)}{n}  \\
&=& \int_1^x t \cdot dR(t) \nonumber \\
&=&xR(x)-R(1)- \int_1^x R(t) dt  \nonumber \\
&=& o(x) \nonumber,
\end{eqnarray}
where $R(t)=\sum_{n \leq t}f(n)n^{-1} $. By the definition of the mean value of a function in (\ref{eq3399.003}), this confirms that $f(n)$ has mean value zero.         
\end{proof}   

This is standard material in the literature, see \cite[ p.\ 4]{PS12}. The second case is covered by a few results on divergent series, which are considered next. \\

Let  $f(n)=\sum_{d \mid n}g(d)$, and let the series $c=\sum_{n \leq x}g(n)n^{-1} $ be absolutely convergent. Under these conditions, the mean value of the function $f(n)$ can be determined indirectly from the properties of the function $g(n)$.

\begin{thm} \label{thm3399.03}    {\normalfont (Wintner)}  Consider the arithmetic functions $f,g : \N \longrightarrow \C$, and assume that the associated generating series are zeta multiple
\begin{equation} \label{eq3399.043}
\sum_{n \geq 1}\frac{f(n)}{n^s} =\zeta(s) \sum_{n \geq 1}\frac{g(n)}{n^s}.
\end{equation}
Then, the followings hold.
\begin{enumerate} [font=\normalfont, label=(\roman*)]
\item If the series $\sum_{n \geq 1}f(n)n^{-s}$ is defined for $\Re e(s) > 1$, and the series $c=\sum_{n \leq x}g(n)n^{-1} $ is absolutely convergent, then the mean value, and the partial sum are given by
$$
M(f)= \lim_{x \to \infty} \frac{1}{x}\sum_{n \geq 1}\frac{g(n)}{n}  \quad \mbox{ and } \quad \sum_{n \leq x}f(n)=cx+o(x).
$$
\item If the series $\sum_{n \geq 1}f(n)n^{-s}$ is defined for $\Re e(s) > 1/2$, and the series $c=\sum_{n \leq x}g(n)n^{-1} $ is absolutely convergent, then, for any $\varepsilon > 0$, the mean value, and the partial sum are given by
$$
M(f)= \lim_{x \to \infty} \frac{1}{x}\sum_{n \geq 1}\frac{g(n)}{n}  \quad \mbox{ and } \quad \sum_{n \leq x}f(n)=cx+O\left (x^{1/2+\varepsilon}\right ).
$$
\end{enumerate}
\end{thm}

\begin{proof} (i) The partial sum of the series (\ref{eq3399.043}) is rearranged as 
\begin{eqnarray}
\sum_{n \leq x}f(n) &=& \sum_{n \leq x} \sum_{d \mid n} g(d)  \\
&=& \sum_{d \leq x} g(d)\sum_{n \leq x/d} 1  \nonumber \\
&=& \sum_{d \leq x} g(d) \left ( \frac{x}{d}- \left \{ \frac{x}{d} \right \} \right )  \nonumber \\
&=&x \sum_{d \leq x} \frac{g(d)}{d}+O \left (\sum_{d \leq x}|g(d)| \right ) \nonumber,
\end{eqnarray}
wheere $\{x\}=x-[x]$ is the fractional part function. The first line arises from the convolution of the power series $\zeta(s)=\sum_{n \geq 1}n^{-s}$, and $\sum_{n \geq 1}g(n)n^{-s}$. This is followed by reversing the order of summation. Moreover, the first finite sum is
\begin{equation} \label{eq3399.048}
\sum_{n \leq x}\frac{g(n)}{n} = \sum_{n \geq 1}\frac{g(n)}{n}+o(1) =c +o(1).
\end{equation}
because $\sum_{n \geq 1}g(n)n^{-1}$ is absolutely convergent. A use a dyadic method to split the second finite sum as
\begin{eqnarray}
\sum_{d \leq x}|g(d)|  &=& \sum_{d \leq x^{1/2}}\frac{|g(d)|}{d} \cdot d+ \sum_{ x^{1/2} \leq d \leq x }\frac{|g(d)|}{d} \cdot d \\
&\leq & x^{1/2}\sum_{d \leq x^{1/2}}\frac{|g(d)|}{d}+ x \sum_{ x^{1/2} \leq d \leq x }\frac{|g(d)|}{d}  \nonumber \\
&=& O\left ( x^{1/2} \right ) +o(x) \nonumber\\
&=& o(x) \nonumber.
\end{eqnarray}
Again, this follows from the absolute convergence $\sum_{n \geq 1}|g(n)|n^{-1}<\infty$.   
\end{proof}

Similar proofs appear in \cite[p.\ 138]{PA88}, \cite[p. 83]{DL12}, and \cite[p.\ 72]{HA05}. Another derivation of the Wintner Theorem from the Wiener-Ikehara Theorem is also given in  \cite[p. \ 139]{PA88}.

\begin{thm} \label{thm3399.07}    {\normalfont (Axer)}  Let $f,g : \N \longrightarrow \C$, be arithmetic functions and assume that the associated generating series are zeta multiple
\begin{equation} \label{eq3399.143}
\sum_{n \geq 1}\frac{f(n)}{n^s} =\zeta(s) \sum_{n \geq 1}\frac{g(n)}{n^s}.
\end{equation}
If the series $\sum_{n \geq 1}f(n)n^{-s}$ is defined for $\Re e(s) > 1$, and $\sum_{n \leq x}|g(n)|=O(x) $ is convergent, then the mean value, and the partial sum are given by
\begin{equation} \label{eq3399.088}
M(f)= \sum_{n \geq 1}\frac{g(n)}{n}  \quad \mbox{ and } \quad \sum_{n \leq x}f(n)=cx+o(x).
\end{equation}
\end{thm}

The goal of the next result is to strengthen Wintner Theorem by removing the absolutely convergence condition.

\begin{thm} \label{thm3399.09}    Let $f,g : \N \longrightarrow \C$, be arithmetic functions and assume that the associated generating series are zeta multiple
\begin{equation} \label{eq3399.243}
\sum_{n \geq 1}\frac{f(n)}{n^s} =\zeta(s) \sum_{n \geq 1}\frac{g(n)}{n^s}.
\end{equation}
If the series $\sum_{n \geq 1}f(n)n^{-s}$ is defined for $\Re e(s) > 1/2$, and the series $c=\sum_{n \leq x}g(n)n^{-1} $ is convergent, then,
\begin{equation} \label{eq3399.088}
M(f)= \sum_{n \geq 1}\frac{g(n)}{n}  \quad \mbox{ and } \quad \sum_{n \leq x}f(n)=cx+o(x).
\end{equation}
\end{thm}

\begin{proof} (i) The partial sum of the series (\ref{eq3399.243}) is rearranged as 
\begin{eqnarray}
\sum_{n \leq x}n \cdot f(n) &=& \sum_{n \leq x}n \sum_{d \mid n} g(d)  \\
&=& \sum_{d \leq x} g(d)\sum_{n \leq x/d} n  \nonumber \\
&=& \sum_{d \leq x} g(d) \left ( \frac{x^2}{2d}+o \left ( \frac{x^2}{2d} \right ) \right )  \nonumber \\
&=&\frac{x^2}{2}  \sum_{d \leq x} \frac{g(d)}{d}+o \left (\frac{x^2}{2} \left |\sum_{d \leq x}\frac{g(d)}{d} \right | \right ) \nonumber,
\end{eqnarray}
wheere $\{x\}=x-[x]$ is the fractional part function. The first line arises from the convolution of the power series $\zeta(s)=\sum_{n \geq 1}n^{-s}$, and $\sum_{n \geq 1}g(n)n^{-s}$. This is followed by reversing the order of summation. Moreover, the first finite sum is
\begin{equation} \label{eq3399.148}
\sum_{n \leq x}\frac{g(n)}{n} = \sum_{n \geq 1}\frac{g(n)}{n}+o(1) =c +o(1).
\end{equation}
because $\sum_{n \geq 1}g(n)n^{-1}$ is absolutely convergent. A use a dyadic method to split the second finite sum as
\begin{eqnarray}
\sum_{d \leq x}|g(d)|  &=& \sum_{d \leq x^{1/2}}\frac{|g(d)|}{d} \cdot d+ \sum_{ x^{1/2} \leq d \leq x }\frac{|g(d)|}{d} \cdot d \\
&\leq & x^{1/2}\sum_{d \leq x^{1/2}}\frac{|g(d)|}{d}+ x \sum_{ x^{1/2} \leq d \leq x }\frac{|g(d)|}{d}  \nonumber \\
&=& O\left ( x^{1/2} \right ) +o(x) \nonumber\\
&=& o(x) \nonumber.
\end{eqnarray}
Again, this follows from the absolute convergence $\sum_{n \geq 1}|g(n)|n^{-1}<\infty$. Lastly, but not least, the original partial sum  is recovered by partial summation.         
\end{proof}

\section{Extension To Arithmetic Progressions} \label{s4902}
The mean value of theorem arithmetic functions over arithmetic progressions $\{ qn + a : n \geq 1 \}$ facilitates another simple proof of Dirichlet Theorem.

\begin{thm} \label{thm4902.31}    {\normalfont }      Let $f(n)=\sum_{d \mid n} g(d)$, and let the series $\sum_{n \geq 1}g(n)n^{-1}\ne0$ be absolutely convergent. Then
\begin{equation} \label{eq257.03}
M(f)= \lim_{x \to \infty}\,\frac{1}{x}\sum_{\substack{n \leq x \\
n \equiv a \bmod q}}f(n)= \frac{1}{q} \sum_{d \mid q} \frac{c_d(a)}{d} \sum_{ n \geq 1} \frac{g(dn)}{n}.
\end{equation}
where $c_k(n)=\sum_{\gcd(x,k)=1}e^{i2 \pi nx/k}$.
\end{thm}

For the parameter $1 \leq a < q$, and $\gcd(a, q) = 1$, the mean value reduces to
\begin{equation} \label{eq257.03}
M(f)= \lim_{x \to \infty}\,\frac{1}{x}\sum_{\substack{n \leq x \\
n \equiv a \bmod q}}f(n)= \frac{1}{q} \sum_{d \mid q} \frac{\mu(d)}{d} \sum_{ n \geq 1} \frac{g(dn)}{n}.
\end{equation}
The proof is given in \cite[p.\ 143]{PA88}, seems to have no limitations on the range of values of $q \geq 1$. Thus, it probably leads to an improvement on the Siegel-Walfisz Theorem, which states that 
\begin{equation} \label{eq257.22}
\pi(x,q,a)=\frac{1}{\varphi(q) }\frac{x}{ \log x}+O \left ( e^{-(\log x)^{\beta}} \right )
\end{equation}
 where $q = O\left ( \log^B x \right )$, with $B > 0$, and $0 < \beta <1$ constants.

\section{Powers Sums Over Arithmetic Progressions}
Let $\N=\{0,1,2,3, \ldots\}$ be the set of nonnegative integers, and let $q |N+a=\{qn+a:n \in \N\}$ be the arithmetic progression defined by a pair of integers $a> 0$, and $q>\geq 1$. The sums of powers over arithmetic progressions is one of the possible generalizations of the sums of powers , $k ? 0$, over the integers. A few estimates of the powers sums over arithmetic progressions are computed here.

\begin{lem} \label{lem66.1} Let $a \geq 0$ and $q\geq 1$ be fixed integers.  Let $x \geq 1$ be a sufficiently large real number. Then
\begin{enumerate} [font=\normalfont, label=(\roman*)]
\item $ \displaystyle \sum_{\substack{n \leq x\\ n \equiv a \bmod q}}n=\frac{1}{2q}x^2+o\left( \frac{1}{q}x^2\right ).$
\item $ \displaystyle \sum_{\substack{n \leq x\\ n \equiv a \bmod q}}n \geq  \frac{1}{2q}x^2+O\left( \frac{1}{q}x\right ).$
\end{enumerate}
\end{lem}

\begin{proof} The integers in a linear arithmetic progression are of the form $n = qm + a$, with $0\leq  m \leq (x- a)/q$. Inserting this into the finite sum produces
\begin{eqnarray} \label{eq257.90}
   \sum_{\substack{n \leq x\\ n \equiv a \bmod q}}n&=& q\sum_{m \leq (x-a)/q)}m +a \sum_{m \leq (x-a)/q)}1 \\
&=& \frac{q}{2} \left [ \frac{x-a}{q} \right ] \left ( \left [ \frac{x-a}{q} \right ] +1 \right ) +a \left [ \frac{x-a}{q} \right ] \nonumber,
\end{eqnarray}
where $[z]=z-\{z\}$ be the largest integer function. Set $z=x-a$, and expand the expression to obtain:
\begin{eqnarray} \label{eq257.92}
&& \frac{q}{2} \left [ \frac{x-a}{q} \right ] \left ( \left [ \frac{x-a}{q} \right ] +1 \right ) +a \left [ \frac{x-a}{q} \right ] \\
&=&\frac{q}{2} \left ( \frac{z}{q}-\left \{ \frac{z}{q} \right \} \right )  \left (  \frac{z}{q}-\left \{ \frac{z}{q} \right \} +1 \right ) +a\left ( \frac{z}{q}-\left \{ \frac{z}{q} \right \} \right )\nonumber \\
&=& \frac{1}{2q}z^2-z\left \{ \frac{z}{q} \right \}+\frac{q}{2}\left \{ \frac{z}{q} \right \}^2+\frac{z}{2}-\frac{q}{2}\left \{ \frac{z}{q} \right \}  +a\left ( \frac{z}{q}-\left \{ \frac{z}{q} \right \} \right ) \nonumber\\ 
&=&\frac{1}{2q}z^2+o\left( \frac{1}{q}z^2\right ) \nonumber.
\end{eqnarray}
Replacing $z=x-a$ back into the (\ref{eq257.92}) yields the result.  
\end{proof}

The above estimates are sufficient for the intended applications. A sharper estimate of the form 
\begin{equation} \label{eq257.33}
 q\sum_{n \leq (x-a)/q}n=\frac{1}{2q}x^2+O\left( \frac{1}{q}x\right )+\frac{1}{4}
\end{equation}
appears in \cite[p.\ 83]{SH83}.

\end{document}